\documentclass[]{interact}

\usepackage{epstopdf}%
\usepackage[caption=false]{subfig}%
\usepackage{listings} %
\usepackage{tikz} %
\usetikzlibrary{positioning,automata}

\usepackage[natbibapa,nodoi]{apacite}%
\theoremstyle{plain}%

\theoremstyle{definition}

\theoremstyle{remark}

\usepackage{color}

\newcommand{\secref}[1]{\S\ref{#1}}
\newcommand{\figref}[1]{Fig.~\ref{#1}}
\newcommand{\ie}{\textit{i.e.}}
\newcommand{\eg}{\textit{e.g.}}
\newcommand{\re}{\operatorname{Re}}
\newcommand{\im}{\operatorname{Im}}

\makeatletter
\def\+{\@postfix+}
\def\*{\@postfix*}
\def\?{\@postfix?}
\def\@postfix#1{{#1}\@ifnextchar){}{\;}}
\makeatother

\begin{document}

\title{A Rewriting System for Convex Optimization Problems}

\author{
\name{
        Akshay Agrawal\textsuperscript{a}\thanks{Akshay Agrawal: 
        akshayka@cs.stanford.edu},
        Robin Verschueren\textsuperscript{b}\thanks{Robin Verschueren: robin.verschueren@imtek.uni-freiburg.de},
        Steven Diamond\textsuperscript{a}\thanks{Steven Diamond: diamond@cs.stanford.edu}, and
        Stephen Boyd\textsuperscript{a}\thanks{Stephen Boyd: boyd@stanford.edu}
    }
\affil{
        \textsuperscript{a}Depts.\ of Computer Science and Electrical
            Engineering, Stanford University, USA; \\
        \textsuperscript{b}Institut f\"{u}r Mikrosystemtechnik,
            Albert-Ludwigs-Universit\"{a}t Freiburg, Germany
    }
}

\maketitle
\begin{abstract}
We describe a modular rewriting system for translating optimization problems
written in a domain-specific language to forms compatible with low-level solver
interfaces. Translation is facilitated by reductions, which accept a category
of problems and transform instances of that category to equivalent instances of
another category. Our system proceeds in two key phases: analysis, in which we
attempt to find a suitable solver for a supplied problem, and canonicalization,
in which we rewrite the problem in the selected solver's standard form.
We implement the described system in version 1.0 of CVXPY, a
domain-specific language for mathematical and especially convex optimization.
By treating reductions as first-class objects, our method makes it easy to
match problems to solvers well-suited for them and to support solvers with a
wide variety of standard forms.
\end{abstract}


\section{Introduction}\label{intro}
Mathematical optimization centers on the optimization problem.
Every optimization problem has three attributes: a
variable whose value is to be assigned,
constraints that the variable must satisfy, and a real-valued objective function
that measures the displeasure or cost incurred by any particular assignment to
the variable. To solve an optimization problem is to find a numerical
assignment to the variable that minimizes the objective function
among all choices that satisfy the constraints.

Unfortunately, most optimization problems cannot be solved
efficiently \citep[\S1.4]{boyd2004}. There are, however, classes of optimization
problems that can be solved in polynomial time. An important such class contains
\textit{convex} optimization problems ---
problems where the objective function is convex and where the constraints
are described by a set of equality constraints with affine functions and 
inequality constraints with convex functions \citep{NN1994, boyd2004}.

Modern convex optimization has its origin in linear programming,
which traces back to the late 1940s,
after the Second World War \citep[\S2]{dantzig1963}. Since then,
convex optimization has been extended to include a much wider 
variety of problems, and has found application in machine learning
\citep{hastie2009}, control \citep{boyd1994}, and computer science
\citep{bertsekas1991, goemans1995, parrilo2003}, to name just a few of the
fields touched by it. To accommodate the applications of convex optimization,
researchers and practitioners have over the years authored many software
packages. We distinguish between two types of software packages:
domain-specific languages, which streamline the process of specifying
optimization problems, and low-level numerical solvers, which furnish solutions
to problem instances.

\subsection{Domain-specific languages}\label{dsl}
A domain-specific language (DSL) is a language that is designed
for a particular application domain \citep{mernik2005}; familiar examples
include MATLAB and SQL. DSLs for convex optimization are languages designed
for specifying convex optimization problems in natural,
human-readable forms, and they obtain solutions to problems on their users'
behalf by invoking numerical solvers; popular ones include Yalmip
\citep{yalmip}, CVX \citep{cvx}, Convex.jl \citep{convexjl}, and CVXPY
\citep{cvxpy}. These DSLs do support some nonconvex
regimes (e.g., combinatorial optimization), and there also exist DSLs
for nonlinear optimization \citep[see][\S1.4, for a definition]{boyd2004},
including GAMS \citep{gams}, AMPL \citep{ampl}, and JuMP \citep{jump};
here, however, we limit our discussion to convexity.

We present below a toy example of an optimization
problem written in CVXPY (version 1.0), a Python-embedded DSL \citep[see][for background on
embedded DSLs]{hudak1996}. The choice of language here is not particularly
important, as the code would look similar if translated to any of the
other aforementioned DSLs.

\begin{verbatim}
    from cvxpy import *

    alice = Variable()
    bob = Variable()

    objective = Minimize(maximum(alice + bob + 2, -alice - bob))
    constraints = [alice <= 0, bob == -0.5]
    toy = Problem(objective, constraints)
    opt = toy.solve()
\end{verbatim}

The CVXPY problem \texttt{toy} has two scalar optimization
variables, \texttt{alice} and \texttt{bob}. Every \texttt{Variable} object
has stored in its \texttt{value} field a numeric value, which is
unspecified upon creation; \texttt{alice} and \texttt{bob} can hold floating
point values. The objective is to minimize a piecewise-affine function of
\texttt{alice} and \texttt{bob}, where the function is represented with the
\texttt{max} \textit{atom}. Atoms are mathematical functions like
\texttt{square} and \texttt{exp} that operate on CVXPY expressions.
CVXPY implements as library functions dozens
of atoms for users to use in constructing problems.
The arguments to the \texttt{max} atom are
\texttt{Expression} objects, which encode mathematical expressions.
\texttt{Constraint} objects are created by linking two expressions with a
relational operator (\texttt{<=}, \texttt{>=}, or \texttt{==}). 
In the second-to-last line, the CVXPY problem \texttt{toy} is constructed,
but not solved.
Finally, an
invocation to \texttt{toy}'s \texttt{solve} method solves the
problem. A side-effect is that the \texttt{value} fields of the optimization
variables present in the problem (\texttt{alice} and \texttt{bob}) are assigned
values that minimize the objective while satisfying the constraints, and the
return value of such a solve is the value of the objective function evaluated
at the variable values. After invoking \texttt{solve} above, we find
that \texttt{alice.value~==~-0.5}, \texttt{bob.value~==~-0.5}, and
\texttt{opt~==~1.0}.  These values satisfy the two constraints, and among all
such assigments, yield the smallest value of the objective function.

The \texttt{solve} method cannot solve all problems: Whether or not a problem
can be solved depends on the objective, constraints, and variables present in
the problem. In particular, these entities must be constructed in such a way
that CVXPY can detect that their assemblage is in fact a convex problem.
Recognizing convexity in general can be difficult.  Many tricks aid in this task, but
their application is sometimes guided by an intuition that is difficult to
codify.  Nonetheless, there do exist sufficient, but not necessary, rulesets
for algorithmically detecting convexity. The ruleset employed by CVXPY, CVX,
Convex.jl, and Yalmip is called \textit{disciplined convex programming} (DCP)
\citep{grant2004}. These DSLs require users to express their problems by means
of DCP, and while not all convex problems adhere to DCP, non-DCP convex
problems can in practice be made DCP-compliant with at most a moderate amount
of human effort and expertise.
DCP-compliant problems can be verified as convex using a simple set 
of rules that are readily automated.

Our CVXPY example does not procedurally describe the method of procuring a
solution. Optimization DSLs are in this sense \textit{declarative languages}.
An optimization problem is simply a precise articulation of preferences and
constraints among the values of the variables;
how the problem is to be solved is another story. This story is
the topic of our next section. But an abridged version is as follows: A
subroutine inspects the problem and invokes a
numerical solver capable of solving it, and most DSLs also allow users to
mandate that a specific solver be used to solve any given problem --- in CVXPY,
for example, users may select a solver via the \texttt{solve} method's keyword
argument \texttt{solver}.

\subsection{Numerical solvers}
A numerical solver is a low-level tool that takes as input an optimization
problem encoded in a rigid format and returns a solution for it. Every solver
is tied to one or more \textit{classes} of problems, insofar as problems
supplied to a solver must be instances of one of its supported classes. One of
the oldest problem classes is least squares, dating back to works authored by
Legendre and Gauss in the late 18th and early 19th centuries (\citealp[see][for
a discussion of Legendre and Gauss' contributions to the methodology]{ls}, and
\citealp[for a translation of Gauss' manuscript on the same]{gauss1995}). Other
well-studied
convex optimization classes include linear programs, popularized by Dantzig
following the Second World War \citep{dantzig1963}, and cone programs,
introduced in the late 20th century by \citeauthor{NN1992} (\citealp{NN1992};
\citealp[\S4.1]{NN1994}).

\begin{figure}
\centering
    \begin{minipage}[t]{.42\textwidth}
        \centering
        \includegraphics[scale=0.40]{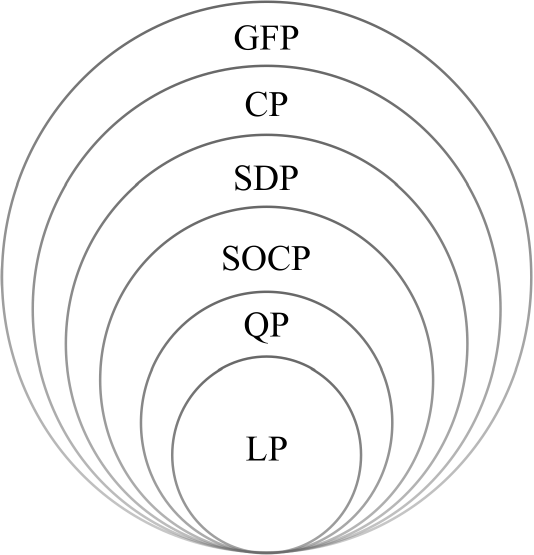}
        \caption{Hierarchy of convex optimization problems. Rewriting systems
        should reduce instances of classes higher in the hierarchy to lower
        ones, when possible. (LP: linear program, QP: quadratic program, SOCP:
        second-order cone program, SDP: semidefinite program, CP: cone
        program, GFP: graph form program.)}
        \label{fig-hierarchy}
    \end{minipage}%
    \hfill
    \begin{minipage}[t]{0.54\textwidth}
    \centering
    \includegraphics[scale=0.40]{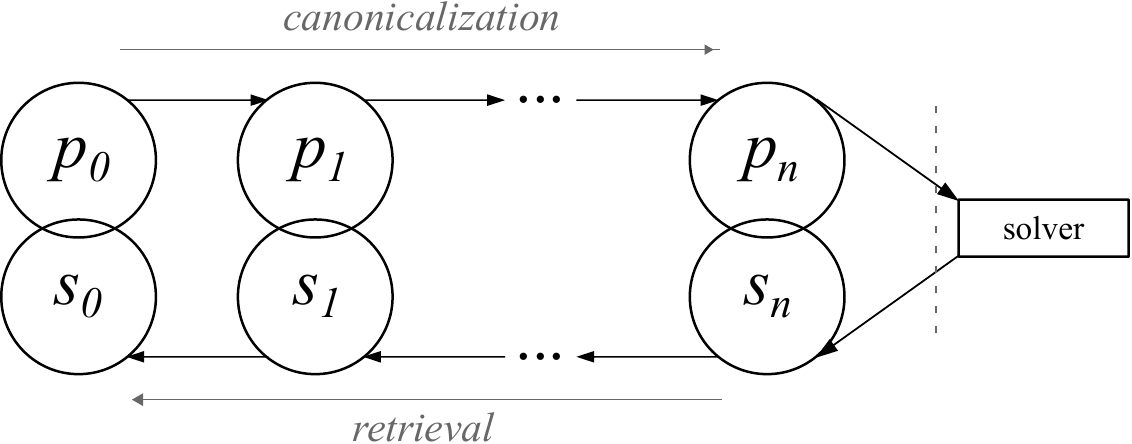}
    \caption{Generic representation of a rewriting,
             capped with an invocation to a solver. The user-posed problem
             $p_0$ is \textit{canonicalized} to a form $p_n$ compatible with
             the targeted solver via a sequence of $n$ reductions, yielding the
             intermediate problems $p_1$ through $p_{n-1}$. The solver solves
             $p_n$ and returns a solution $s_n$ for $p_n$; we then
             \textit{retrieve} solutions for the intermediate problems in
             reverse order, terminating with a solution $s_0$ for $p_0$.}
    \label{fig-canon}
    \end{minipage}
\end{figure}

Various classes of convex problems fit into a hierarchy, as
depicted in \figref{fig-hierarchy}. Every linear program reduces to a quadratic
program \citep[\S4.4]{boyd2004}; every quadratic program reduces to a
second-order cone program \citep[\S6.2.3]{NN1994}; every second-order cone
program reduces to a semidefinite program \citep{sdp}; every
semidefinite program reduces to a cone program; and every cone program reduces
to a graph form program \citep{parikh2014}. Problem classes increase in both
generality and difficulty as one goes up the hierarchy. For example, whereas
every second-order cone program is also a semidefinite program, it is not
advisable to solve the second-order cone programs using semidefinite program
solvers \citep{alizadeh2001}. This fact motivates the existence of numerical
solvers for different classes of problems --- all other things equal, it is
better to use as specific a solver as possible for the problem at hand.
 
In 1953, Hays and Dantzig developed for the RAND Corporation one of
the earliest linear program solvers \citep[\S2.1]{dantzig1963}. The universe
of solvers has since grown in lockstep with the universe of optimization
algorithms. There exist today solvers for many types of convex programs.
Examples of solvers specialized to each class include GLPK \citep{GLPK} and CBC
\citep{cbc} for linear programs;
qpOASES \citep{qpoases} and OSQP \citep{osqp} for quadratic programs;
ECOS \citep{ecos} and Gurobi \citep{gurobi} for second-order cone
programs; MOSEK \citep{mosek} and SeDuMi \citep{Sturm1999} for semidefinite programs;
SCS \citep{scs} for more general cone programs;
and POGS for graph form programs \citep{pogs}.
Many of the solvers listed also support nonconvex problems such as mixed-integer
and nonlinear programs.

When we say that numerical solvers are low-level tools, we mean that it is
onerous to translate problems to forms acceptable to solvers --- even
deciding which solver to use for a particular problem is a skill that requires
training. These observations provide two of the raisons d'\^{e}tre for
optimization DSLs.

\subsection{Canonicalization}
The process of converting an optimization problem encoded in a DSL to a
solver-compatible form --- for example, the process by which CVXPY transformed
our toy problem to the above representation --- is called
\textit{canonicalization} \citep[\S4]{grant2004}. Consider once again the toy
problem from \secref{dsl}, the one transcribed in CVXPY. This
problem can be transformed to an equivalent linear program, \ie, a problem of
the form \begin{equation*}
\begin{array}{ll}
\mbox{minimize} & c^Tx \\
\mbox{subject to} & Gx \leq h \\
& Ax = b,
\end{array}
\end{equation*}
where $x$ is the (vector) variable, the matrices $G$ and $A$, and vectors
$c$, $h$, and $b$, are constants, and
the inequality is component-wise. Upon invoking the solve method,
CVXPY canonicalizes the toy problem to the above standard form to make it
compatible with a numerical (linear program) solver.
The transformed problem has a variable $x
\in \mathbf{R}^3$ whose first and second components represent \texttt{alice}
and \texttt{bob}, respectively, and whose third component is an auxiliary
variable introduced in the canonicalization. The problem data are
\begin{align*}
G = \begin{bmatrix}
1 & 1 & -1 \\
-1 & -1 & -1 \\ 
1 & 0 & 0 
\end{bmatrix}, \quad
A = \begin{bmatrix} 0 & 1 & 0 \end{bmatrix}, \quad
c = \begin{bmatrix} 0 \\ 0 \\ 1 \end{bmatrix}, \quad
h = \begin{bmatrix} -2 \\ 0 \\ 0 \end{bmatrix}, \quad
b = -0.5.
\end{align*}

This canonical form of our toy problem is the result of applying the
transformation described in \S6.5 of \citeauthor{dantzig1997}'s
\citeyear{dantzig1997} text, and as such it could have easily been produced
manually. Yet as problems of interest grow larger and more complex, producing
canonical forms by hand quickly becomes a tedious, laborious, and error-prone
task. Many users of convex optimization still canonicalize problems manually
--- instead of letting optimization DSLs do the work for them --- due to
reliance on legacy systems, performance concerns, or in some cases,
ignorance of the existence of DSLs that can automate the task.

Problems written in CVX, Convex.jl, and CVXPY are automatically canonicalized
to conic form, and Yalmip supports other classes of problems as well. For all
of these DSLs, however, the solver selection and canonicalization procedures
are implemented in ad hoc fashions that cannot easily be modified or extended,
for instance to target new problem classes.

\subsection{This paper}\label{this-paper}
This paper is about treating canonicalization as a first-class component of
the software ecosystem for convex optimization. We propose a
\textit{rewriting system} that sits between DSLs
and numerical solvers, translating problems expressed in the former to
forms compatible with the latter.

The atomic rewriting unit in our rewriting system is the \textit{reduction}.
 A reduction is a function that converts problems of one form to
equivalent problems of another form. Two problems are \textit{equivalent} if a
solution for one can be readily converted to a solution for the other \citep[as
defined in][\S4.1.3]{boyd2004}; readers from computer science will
recognize that our definitions are in the same spirit as the more formal ones
from their field (\citealp[see for example \S10.3 of ][]{papadimitriou1994},
and \citealp[\S2.2 of][]{arora09}).

Every canonicalization is a reduction. In many cases, canonicalizations are
complex enough to merit decomposing them into compositions of reductions. For a
rewriting system to be useful, it must \textit{retrieve} solutions for
canonicalized problems to solutions for their provenances; because reductions
output equivalent problems, they by definition support retrieval.
The dual processes of canonicalization and retrieval are diagrammed in
\figref{fig-canon}.  Any rewriting that is cast as a composition of reductions
is provably correct: The rewriting will either (a) output an equivalent
problem, if the reductions
in the canonicalization are mutually compatible and the first reduction is
applicable to the source problem, or otherwise (b) audibly fail. There are at
least three other benefits enjoyed by placing special emphasis on problem
rewritings and reductions: doing so provides a structured way of preferentially
targeting some solvers over others,
simplifies the interfacing of
domain-specific languages with new solvers, and unifies problem transformations
and back-end optimizations like presolves within a single conceptual
framework.

The remainder of this paper is structured as follows. In
\secref{architecture}, we elaborate upon the role and structure of rewriting
systems; in \secref{examples}, we list many examples of
reductions, from simple to nuanced; and in \secref{impl}, we discuss version
1.0 of CVXPY, an open-source implementation of our proposed rewriting system.

\section{An architecture for rewriting systems}\label{architecture}
\subsection{Principles}\label{design}
There are four principles to which all optimization rewriting systems should
adhere; these principles are informed by ones adopted by the software
compiler community \citep[\S1.4.2]{dragonbook}.

\begin{enumerate}
    \item Every rewriting must yield an equivalent problem that is target-compatible.
    \item For each problem, an effort should be made to select a
          suitable solver for it.
    \item The rewriting time must be tolerable.
    \item The engineering effort required to maintain the rewriting system
          and add solvers to it must be kept manageable.
\end{enumerate}
Adherence to the first three principles is necessary in order for a rewriting
system to be useful, and to faithfully solve the specified problem.
The fourth principle expedites not only the work of
engineers responsible for rewriting systems but also that of those developing
new solvers: Researchers who interface their solvers to popular DSLs gain
immediate access to a rich ecosystem of problems for testing and tuning
their algorithms. Indeed, the advent of software compilers fundamentally
altered the development cycle of processors
--- today, compilers are built before processor designs are finalized, for
often a processor is only useful if a compiler can exploit it
\citep[\S1.5.3]{dragonbook}. The creation of optimization rewriting systems
that satisfy the fourth principle might effect a similar paradigm shift in the
design of numerical solvers.

As for how to actually satisfy these principles --- the first principle is
automatically satisfied by the use of reductions. One way to satisfy the
remaining three principles is to decompose the rewriting system architecture
into three distinct phases.  

\subsection{Phases}\label{phases}
\begin{figure}
\centering
\includegraphics[width=\textwidth]{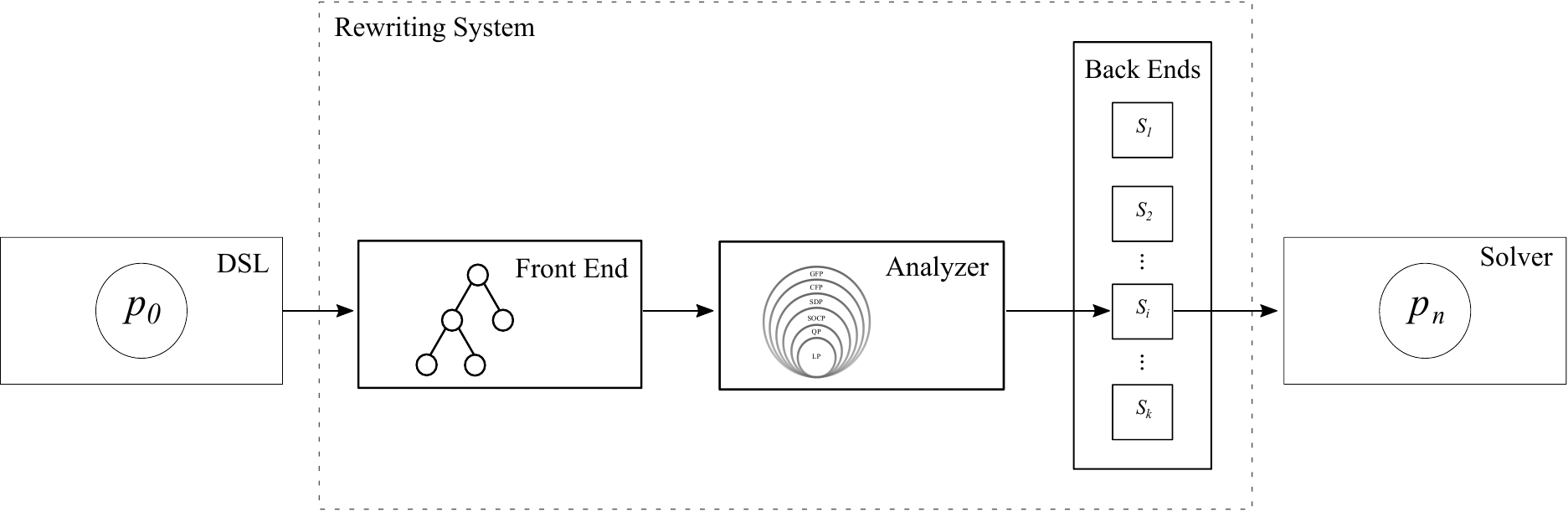}
\caption{Architecture of a three-phase rewriting system. The system
takes as input a problem $p_0$ written in a DSL\@. A \textit{front end} specific
to that DSL takes $p_0$ and encodes it in some intermediate representation,
e.g., mathematical expression trees. The \textit{analyzer} reads the
intermediate representation, discovers what it can about the problem's
structure, and then selects a target, or solver;
the analyzer may itself produce further intermediate representations
of the problem via reductions. A \textit{back end} takes the
intermediate representation produced by the analyzer and reduces it to a form
$p_n$ that is compatible with the target.
}
\label{fig-phases}
\end{figure}
Just as it is natural to decompose software compilers into three phases
\citep[\S1.2]{dragonbook}, so it is with optimization rewriting systems
(\figref{fig-phases}). The first and third phases of both systems are
analogous.  In the first phase, a \textit{front end} takes a human-readable
specification of a program and converts it to an intermediate representation;
both compilers and rewriting systems often use abstract syntax trees as their
intermediate representations (\citealp[\S2.5.1]{dragonbook};
\citealp{convexjl}, \citealp{DB:15}). In the third phase, a \textit{back end} takes an
intermediate representation and translates it to a target-compatible form;
targets for software compilers are processor architectures, while targets for
optimization rewriting systems are numerical solvers. Our front end does not
make use of reductions, as the problem of constructing a syntax tree is one of
parsing, but our back ends most certainly do --- in fact, our back ends are
nothing more than sequences of reductions.
In other words, each of the back ends supported by the rewriting system is
a different canonicalization procedure.

Whereas the second phase of a compiler optimizes programs
\citep[\S1.2.5]{dragonbook}, the second phase of our rewriting system
\textit{analyzes} problems, selecting for each a suitable target; an imperfect
analogy might compare our analysis phase to a static analysis tool like
Frama-C \citep{frama-c}. Compilers like the GNU Compiler Collection \citep{gcc}
can translate a given source program to any of their targeted architectures,
since every architecture is Turing complete \citep{sipser1996}. Rewriting
systems provide no such guarantees --- solvers are specific to particular
classes of problem, \ie, a given optimization problem might be compatible
with but a subset of the solvers targeted by a rewriting system, hence the need
for analysis. A sensible analysis policy is to identify the most specific class
to which a problem belongs and choose a target supporting that class (see
\figref{fig-hierarchy}).  Solver selection might also reflect desiderata like
accuracy, scalability, and speed. For example, first-order solvers like SCS
scale to larger problems than interior-point methods like ECOS, though the
latter typically provide more accurate solutions than do the former.
Analysis may of course fail to find a target that can handle a particular
problem.  When this happens, the rewriting system should abort with a
descriptive error code or message.

This three-phase architecture satisfies principles two through four
listed in \secref{design}. The existence of an analysis phase satisfies the
second principle, a best-effort implementation of the analysis phase will
satisfy the third principle, and the separation of back ends from front ends,
together with the use of reductions as modular rewriting units, satisfies the
fourth principle.

\section{Reductions}\label{examples}
In this section, we list several examples of reductions. Some of the listed
reductions might be used to perform routine operations common among
canonicalizations, while other more involved ones might be used to reduce
problems to instances of simpler classes. \\

\noindent
\textit{Notation.} In all of the following examples, the variable $x$ denotes
the optimization variable wherever it appears, which may be a scalar or a
vector. If we refer to the variables, plural, of an optimization problem, we
mean to refer to the individual components of $x$. An
\textit{equality constraint} is one of the form $f(x) = g(x)$. An
\textit{inequality constraint} is one of the form $f(x) \leq g(x)$, where $f$ and
$g$ are (real) vector-valued functions and the inequality
is component-wise; the \textit{constraint function} of a constraint is $f - g$.
Equality constraints are convex if the constraint function is affine, and
inequality constraints are convex if the function on the left ($f$)
is convex and the function on the right ($g$) is concave.

\subsection{Simple reductions}\label{simple} In this section, we present some
simple but useful reductions. \\

\noindent
\textit{Flipping objectives.} The action of
inverting a maximization objective to a minimization is a reduction: Maximizing
a function $f$ over its domain is equivalent to minimizing $-f$, as solutions
to the two problems are related by an inversion of sign. No work is required to
retrieve a solution from the reduced-to problem, as the optimal sets
of both problems are the same.  (The optimal value of the problem must be negated.)\\

\noindent
\textit{Moving expressions to the left of a relation.} Subtracting the
right-hand side of each constraint of a problem from both sides yields
an equivalent problem in which all constraints have zero-valued right-hand
sides. As in the case of the previous reduction, solution retrieval requires no
work. \\

\noindent
\textit{Eliminating general linear inequality constraints.}  Linear inequality 
constraints (\ie, those with $f$ and $g$ both affine) can be replaced with equality
constraints by introducing nonnegative variables. To
wit, the constraint $f(x) \leq g(x)$ holds if and only if $f(x) + s = g(x)$,
for some $s \geq 0$; the auxiliary variable $s$ is called a \textit{slack
variable} \citep[see][\S4.3]{dantzig1963}. Applying this reduction gives an
equivalent problem in which every linear constraint is either an 
equality constraint
or a nonnegativity constraint on a slack variable.
A solution for the original problem may be retrieved from a solution
for the transformed problem by simply fetching the value of $x$ and
discarding the slack variables. \\

\noindent
\textit{Monotone transformations of objective and constraints.}
Composing any monotone increasing function with the objective function of
a problem yields an equivalent problem; so does transforming
any number of constraints by applying any monotone increasing function to both
sides. The retrieval method for this reduction is essentially a no-op, as the
feasible and optimal sets for the two problems are identical.  This reduction
has been employed for centuries --- squaring the Euclidean norm when it appears
as an objective function to render it differentiable is, at least historically,
standard mathematical practice.

Stated in the opposite direction, if the objective function of a convex problem
is monotone increasing, then the problem given by replacing the objective
function with a composition of its inverse and itself (and similarly for the
constraints) is an equivalent convex problem. This reduction
might be used to eliminate exponentials and logarithms from a problem, which in
turn might make the problem amenable to more solvers. \\

\noindent
\textit{Changing variables.} 
Let $\phi$ be any one-to-one function whose image covers the problem domain.
Replacing the optimization variable $x$ wherever it appears with $\phi(z)$ is
a reduction, yielding an equivalent problem with optimization variable $z$.
To retrieve a solution for the original problem from one for the transformed
problem, simply take $x^\star = \phi(z^\star)$, where $x^\star$ denotes an
optimal point for the original problem and $z^\star$ an optimal point for the
transformed problem. 

Changing variables can convert nonconvex problems to equivalent convex
problems. A classic example is the technique for convexifying geometric
programs; this technique both changes variables and transforms the objective
and constraints. A geometric program is an optimization problem of the
form
\begin{equation*}
\begin{array}{ll}
\mbox{minimize}   & f_0(x) \\
\mbox{subject to} & f_i(x) \leq 1, \quad i=1, \ldots, m \\
                  & h_i(x) = 1, \quad i=1, \ldots, p,
\end{array}
\end{equation*}
where the functions $f_i$ are \textit{posynomials} and the functions $h_i$ are
\textit{monomials}.  A monomial is a function $f : \mathbf{R}^n \rightarrow
\mathbf{R}$ over the
nonnegative orthant defined as
\[
f(x) = c x_1^{a_1} x_2^{a_2} \cdots x_n^{a_n},
\]
where $c > 0$ and $a_i \in \mathbf{R}$, and a posynomial is a sum of
monomials. Performing a change of variables with $x_i = \exp(z_i)$ and
taking logarithms of both the objective and constraints results
in a convex problem. In fact, if the $f_i$ are all monomials, then the
resulting problem is a linear program. For a brief introduction to geometric
programming, see \citep[\S4.5]{boyd2004}; for a longer survey, see \citep{geo}.
\\

\noindent
\textit{Eliminating complex numbers.}
It is possible to reduce an optimization problem over a complex
domain to one over a real domain.  (While such a problem has complex
variables and expressions, the constraint and objective functions 
must all be real-valued.)

We provide here a partial specification of the reduction. Absolute values of
complex numbers reduce to Euclidean norms of their concatenated real and
imaginary parts, sums of complex numbers reduce to sums of their real and
imaginary parts, and equality constraints between two complex numbers reduce to
equality constraints on the implicated expressions' real and imaginary parts.
Perhaps more interesting, positive semidefinite constraints on Hermitian
matrices reduce to positive semidefinite constraints on real symmetric
matrices. As presented in \citep{goemans2004}, a Hermitian matrix $X$ is
positive semidefinite if and only if the real symmetric matrix $T(X)$ is
positive semidefinite, where the mapping $T$ is defined as
\[
T(X) = \begin{bmatrix}
            \re X & -\im X \\
            \im X & \re X
       \end{bmatrix}.
\]

As this reduction expands the optimization variable into its real and imaginary
parts, retrieving a solution for the complex-domain problem from a solution for
the real-domain problem requires but a bit of book-keeping to map variable
values from the latter solution to the real and imaginary parts of the original
variable.

\subsection{Presolves}
A presolve is any reduction that is meant to decrease the computational cost 
incurred when solving a problem. Presolves are typically performed
immediately before problems are solved, with some but not all numerical solvers
subjecting problems to a battery of presolves prior to solving them. Many
presolves are applicable across solvers, \ie, a presolve that helps one solver
is likely to help many others. This motivates folding presolves into rewriting
systems and excising them from numerical solvers whenever possible. As
recommended by \citet{bradley1983}, one might even treat the application of
presolves as a fixed point iteration, cyclically applying presolves
until the problem cannot be further simplified. This approach resembles
the multiple passes an optimizing software compiler may make over an
intermediate code representation \citep[\S8 and \S9]{dragonbook}.

There is substantial literature on presolves. \citet{andersen1995}
cast presolves as reductions, listing many examples for linear
programs together with methods for retrieving their solutions. Earlier surveys of
linear programming presolves include \citep{brearley1975}, \citep{bradley1983},
and \citep{tomlin1986}. \citet{tomlin1975} discusses the problem of scaling
data matrices to coax faster performance out of the simplex algorithm, while
\citet{bradley2010} provides a more modern perspective on scaling for a wider
class of algorithms. Here, we present a sample of some of the
presolves covered by these and other references. \\

\noindent
\textit{Eliminating fixed variables.} 
Any variable that is constrained to be a constant is called a fixed variable;
replacing every occurrence of it with the value of the constant yields an
equivalent problem. In the software compilers literature, this technique is
called \textit{constant propagation} \citep[\S9.4]{dragonbook}. Solution
retrieval simply requires setting the values of the fixed variables to their
respective constant values and copying all other variable values.\\

\noindent
\textit{Eliminating free variables.} 
Any variable that does not have upper and lower bounds is called a free
variable; replacing every occurrence of it with the difference of two auxiliary
nonnegative variables yields an equivalent problem for which a solution can be
retrieved in the obvious way \citep[\S 4.3]{dantzig1963}. In the setting of
cone programs, free variables are those that are not restricted to lie in a
cone (other than $\mathbf{R}^n$). There are a number of ways to exploit free
variables in cone programs, some of which are outlined by \citet{anjos2008}. \\

\noindent
\textit{Eliminating redundant constraints.}
Any constraint whose removal leaves the feasible region unchanged is redundant;
deleting such constraints yields an equivalent problem. For example,
any equation in a linear system that is a linear combination of the others
is redundant \citep[\S B.2]{dantzig1997}. As another example, if it is required
that $x \leq b$ and $x \leq c$, and moreover if it is known that
$b \leq c$, then the constraint $x \leq c$ is redundant.
Solution retrieval for this reduction is a no-op.  \\

\noindent
\textit{Scaling.}
Scaling both sides of a constraint by a positive constant is a presolve; this
presolve is an instance of monotonically transforming constraints (see
\secref{simple}).  It has long been known that scaling matrices (\ie, scaling
constraints of the form $Ax \leq b$ or $Ax = b$) can by lowering the
condition numbers of these matrices dramatically affect the performance of
first-order methods for convex optimization (see, \textit{e.g.},
\citealp{tomlin1975}). One scaling technique, called diagonal preconditioning,
premultiplies such constraints by diagonal matrices and also changes variables
by premultiplying the optimization variable by another diagonal matrix
(\citealp[\S2.5]{kelley1995}; \citealp{pock2011}; \citealp{takapoui2016}). A
popular instantiation of this technique is matrix equilibration, which
chooses the diagonal matrices so that all rows of the scaled data matrix have
one $p$-norm and all columns have another, with the two equal for square
matrices.  The literature on equilibration spans decades --- see, for example,
\citep{sluis1969}, \citep{bradley2010}, \citep{mf_equil} and the references
therein.

\subsection{Conic canonicalization of DCP programs}\label{graph-impl}
The embedded languages CVX, Convex.jl, CVXPY, and YALMIP canonicalize problems
to a form compatible with cone program solvers; in particular, the
canonicalized objective function is affine and all constraints are conic,
imposed only on affine expressions of the optimization variable.  These tools
canonicalize problems in the same fashion, and the methodology shared among
them is the subject of this section. The methodology --- which is a reduction
if the problem operated upon is DCP-compliant --- proceeds in three steps:
first, the problem is lifted into a higher dimension via its \textit{Smith
form}, making affine the arguments of each atom; second, the lifted problem is
\textit{relaxed} to a convex problem; and third, every nonlinear atom is
replaced with conic constraints that encode its \textit{graph implementation}.
Our exposition in this section draws from work by \citet{smith1996}, who
introduced Smith form, \citet{grant2008}, who introduced graph implementations,
and \citet{chu2013}, who illustrated these three steps with a clear example. \\

\noindent
\textit{Smith form.}
It is natural to view an optimization problem as composed of a forest of
mathematical expression trees, with one tree for the objective
function and two trees for every constraint, one for each side of the
constraint. The inner nodes of an expression tree represent mathematical
functions, or atoms, and the leaves represent variables and constants. Every
inner node is evaluated at its children, \ie, the children of an atom are its
arguments. For example, the expression $f(x) + c$, in which $f$ is an atom and
$c$ a constant, parses into a tree where the summation atom is the root, $f$
and $c$ are the children of the root, and $x$ is the child of $f$.

Converting an optimization problem to Smith form requires making a single pass
over every expression tree present in the problem. Recursively, beginning with
the root, each subexpression is replaced with an auxiliary variable, and
equality constraints are introduced between the auxiliary variables and the
subexpressions they replaced. The resulting problem is said to be in Smith
form, a key property of which is that the arguments of each atom within the
problem are affine (indeed, they are unadorned variables). Transforming a
problem to Smith form is always a reduction. This reduction does not however
preserve convexity, as any convex function present in the original problem will
appear as the constraint function of an equality constraint in the transformed
problem. \\

\noindent
\textit{Relaxed Smith form.}
If the original problem is DCP-compliant, then its Smith form can be relaxed
to an equivalent convex problem. In particular, relaxing in the correct
direction the nonconvex equality constraints between the auxiliary variables
and their associated atoms is in this case a reduction, in the following sense: if
$(x^\star, t^\star)$ is optimal for the relaxed problem, $t$ a vector of the
auxiliary variables and $x$ the original variable, then $x^\star$ is optimal
for the original problem. \\

\noindent
\textit{Graph implementations.}
The final step in this reduction replaces every constraint in which a
nonlinear convex atom appears with conic constraints that encode the atom's
epigraph.  The \textit{epigraph} of a function is defined as the set of points
that lie above its graph: For a function $f$, its epigraph is defined as the
set of points $\{(x, y) \mid f(x) \leq y\}$. As a simple example of such a
replacement, the epigraph of the function $|x|$ is the set 
$\{(x, y) \mid x \leq y, -x \leq y\}$; accordingly, the constraint $|t_1| \leq
t_2$ would be replaced with the constraints $t_1 \leq t_2$ and $-t_1 \leq
t_2$. Constraints encoding the epigraph of an atom are called the
\textit{graph implementation} of the atom, coined in \citep{grant2004},
though such constraints might more aptly be referred to as an epigraph
implementation; the action of replacing a nonlinear atom with its graph
implementation is called a \textit{graph expansion}.

Graph implementations are useful outside of conic canonicalization as well.
For example, one might choose to only perform graph expansions for
piecewise-linear atoms such as $\texttt{abs}$, $\texttt{max}$, and
$\texttt{sum\_k\_largest}$, which sums the $k$ largest entries of a vector or
matrix. This process of $\textit{eliminating piecewise-linear}$ atoms is
itself a reduction if the problem to which it is applied is DCP-compliant.

\subsection{Other reductions}
The reductions presented are somewhat subtle and problem specific,
examples of the kinds of experimental reductions one might include in a
rewriting system. \\

\noindent
\textit{Decomposing second-order cone constraints.} A \textit{second-order cone
constraint} on a block vector $(x, t)$ is a constraint of the form $\|x\|_2
\leq t$, where $x \in \mathbf{R}^n$ is a vector and $t$ is a scalar. The
dimension of such a second-order cone constraint is $n+1$. Any second-order
constraint of dimension $n+1$, $n \geq 2$, can be reduced to $n-1$
three-dimensional second-order cone constraints by the following chain of
observations: The constraint $\|x\|_2 \leq t$ holds if and only if $x_1^2 +
x_2^2 + \cdots + x_n^2 \leq t^2$ and $t \geq 0$,
which in turn holds if and only if $x_2^2 + \cdots + x_n^2 \leq u^2$,
$u^2 \leq t^2 - x_1^2$, and $u,t \geq 0$, where $u$ is a scalar
variable, \ie, if and only if $(x_2, \ldots, x_n, u)$ and $(x_1, u, t)$ are in
the second-order cone. The result follows by recursing on $(x_2, \ldots, x_n,
u)$. \\

\noindent
\textit{Decomposing semidefinite constraints.} A \textit{semidefinite
program} is a convex optimization problem of the form
\begin{equation*}
\begin{array}{ll}
\mbox{minimize} & \mathbf{tr}(CX) \\
\mbox{subject to} & \mathbf{tr}(A_iX) = b_i, \quad i=1, \ldots, p \\
& X \in \mathbf{S}^n_{+},
\end{array}
\end{equation*}
where the constraint $X \in \mathbf{S}^n_{+}$ requires $X$ to be an $n\times n$ positive
semidefinite matrix. Semidefinite programs can become significantly harder to
solve as the size of the matrices involved increases, motivating reductions
that decompose the matrices into smaller ones.

A semidefinite program that exhibits a chordal aggregate sparsity pattern
can often be reduced to a program involving smaller matrices. The sparsity pattern of
an $n\times n$ matrix is a set $E$ of pairs $(i, j)$, $i, j \in \{1, 2, \ldots,
n\}$, such that $A_{ij} = 0$ if $(i, j) \not\in E$ and $i \neq j$;
the aggregate sparsity pattern of a semidefinite program is the union of the edge sets of the sparsity
patterns of $C$ and $A_1, A_2, \ldots, A_p$. A sparsity pattern is
\textit{chordal} if the induced graph (with vertices $V=\{1,2,\ldots,n\}$ and edges $E$) is chordal,
\ie, if every cycle of length greater than three contains an edge between two non-consecutive
vertices, called a chord.

The reduction in question replaces the optimization variable with smaller
matrix variables, one for each clique in the graph induced by the aggregate sparsity pattern,
rewriting the objective and equality constraints in terms of these new
variables and adding equality constraints to preserve the semantics of the
original problem. The particulars of this reduction (and related ones) can be
found in \citet{fukuda2001}, whose authors were among  the first to exploit
chordal sparsity patterns to decompose large semidefinite programs, and
\citep[\S14.2]{vandenberghe2015}, which thoroughly surveys the topic of chordal
graphs as they relate to semidefinite optimization. \\

\noindent
\textit{Relaxing convex equality constraints.} Consider an optimization problem
with a convex objective, convex inequality constraints, and a single convex
equality constraint. This problem is not convex; however, it can in certain
cases be coerced into a convex form. Letting $x$ be the problem variable, if
there is an index $r$ such that the objective is monotonically increasing in
$x_r$, each inequality constraint function is nondecreasing in $x_r$, and
the equality constraint function is monotonically decreasing in $x_r$, then
relaxing the equality constraint to a nonpositive inequality constraint
produces an equivalent convex problem, \ie, the relaxation is tight
\citep[exercise 4.6]{boyd2004}. 

\subsection{An example}\label{quadratic}
We address in this section the following quite practical question: What types
of convex programs reduce to quadratic programs?

A \textit{quadratic program} is an optimization problem in which the objective
function is a convex quadratic and the constraint functions are affine
\citep[\S4.4]{boyd2004}; quadratic programs have been studied since the 1950s. 
Evidently, a problem in which every inequality
constraint function is piecewise-linear and every equality constraint function
is affine can be reduced to a problem in which every constraint function
is affine (by eliminating the piecewise-linear atoms via graph expansions,
as described in \secref{graph-impl}).

Describing acceptable objective functions requires more work; we will specify
acceptable objective functions via their expression trees, and we will specify
an acceptable expression tree by providing \textit{regular expressions}
\citep[\S3.3.3]{dragonbook} for paths beginning at the root and terminating at
(the ancestor of) a leaf. Letting $A$ denote an affine atom, $P$ a
piecewise-linear atom, and $Q$ a quadratic atom, it is clear that any objective
function whose root-to-leaf paths satisfy the regular expression $A\*QA\*|P\+$
can be canonicalized to a quadratic by eliminating the piecewise-linear
atoms --- this is evident because the Hessian of such a function is constant
almost everywhere.

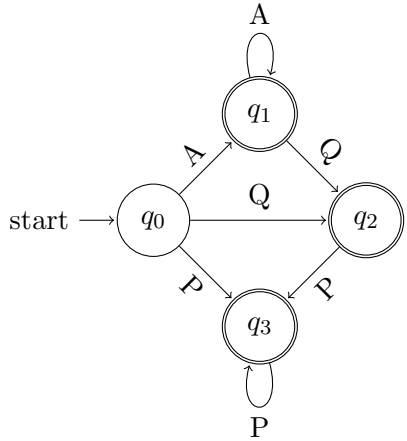
\begin{figure}
\centering
\begin{tikzpicture}[shorten >=1pt,node distance=2cm,on grid]
  \node[state,initial]   (q_0)                {$q_0$};
  \node[state, accepting]           (q_1) [above right=of q_0] {$q_1$};
  \node[state,accepting] (q_2) [below right=of q_1] {$q_2$};
  \node[state,accepting] (q_3) [below right=of q_0] {$q_3$};
  \path[->] (q_0) edge                node [sloped, above] {A} (q_1)
                  edge                node [above] {Q} (q_2)
                  edge                node [sloped, below] {P} (q_3)
            (q_1) edge [loop above]   node         {A} ()
                  edge                node [sloped, above] {Q} (q_2)
            (q_2) edge                node [sloped, below] {P} (q_3)
            (q_3) edge [loop below]   node         {P} ();
\end{tikzpicture}
\caption{A finite-state machine for the example in \secref{quadratic}. Any
DCP-compliant problem in which the root-to-leaf paths of the objective
function's expression tree are accepted by this state machine can be reduced to
a quadratic program, provided that the equality constraint functions are affine
and the inequality constraint functions piecewise-linear. Above, $A$ represents
affine atoms, $Q$ represents quadratic atoms, and $P$ represents
piecewise-linear atoms.} \label{fig-qp-objective}
\end{figure}

Barring nonlinear transformations (\eg, squaring a norm),
one might reasonably assume that the class of problems reducible to quadratic
problem cannot be further generalized, for the regular expression does after all
capture both linear programs and quadratic programs. But such an assumption would
be incorrect. In fact, any DCP-compliant problem (with acceptable constraints)
whose objective function's root-to-leaf paths satisfy the regular expression
$A\*QP\*|P\+$ can be reduced to a quadratic program (\figref{fig-qp-objective}
renders the regular expression as a finite-state machine).  The corresponding
reduction eliminates the piecewise-linear atoms using graph expansions and then
massages the objective into a quadratic.

As an example, the following problem can be canonicalized to a quadratic
program:
\begin{equation*}
\begin{array}{ll}
\mbox{minimize} & (\max(x, 0) + \max(x-1, 0))^2.
\end{array}
\end{equation*}
Note that unlike the simpler class of problems we described, this problem
does not have a constant Hessian --- its second derivative is $0$ when $x < 0$,
$2$ when $0 < x < 1$, and $8$ when $x > 1$. It is nonetheless reducible to a
quadratic program, \ie, a problem whose objective function has a constant
Hessian.

\section{Implementation}\label{impl}
We have implemented a number of the reductions from \secref{examples} in
version 1.0 of CVXPY, an open-source implementation of our proposed three-phase
rewriting system, available at
\begin{center}
    \url{http://www.cvxpy.org/}.
\end{center}

All problem rewriting is facilitated by \texttt{Reduction} objects, and every
reduction implements three methods: \texttt{accepts}, \texttt{apply},
and \texttt{retrieve}. The \texttt{accepts} method takes as input a problem
and returns a boolean indicating whether or not the reduction can be applied
to the problem, the \texttt{apply} method takes as input a problem and
returns a new equivalent problem, and the \texttt{retrieve} method takes
a solution for the problem returned by an invocation of \texttt{apply} and
retrieves from it a solution for its problem of provenance. Some of the tasks
carried out by reductions in our system include flipping objectives
(\secref{simple}), eliminating piecewise-linear atoms (\secref{graph-impl}),
canonicalizing problems to cone programs, (\secref{graph-impl}), and
canonicalizing problems to quadratic programs (\secref{quadratic}). Sequences
of reductions are represented by \texttt{Chain} objects, which are themselves
reductions, and solver back ends are implemented with \texttt{Chain} objects.

Creating expressions and constraints in CVXPY invokes behind-the-scenes
a front end that parses them into expression trees; this functionality
is not new \citep[see][]{cvxpy}. What is new is the method by which solvers are
chosen for problems and the methods by which problems are canonicalized to
their standard forms. In CVXPY 1.0, invoking the \texttt{solve} method of a
problem triggers an analyzer, phase two of our rewriting system. The analyzer
determines the most specific class to which the problem belongs by checking
which back ends accept the problem; back ends are checked in order of
decreasing specificity, and the analysis is short-circuited as soon as
a suitable back end is found. Analysis may itself apply simple reductions
to the problem. Any such reductions are prepended to the back end to create
a \texttt{Chain} object that encapsulates the entire rewriting process, and
the problem is solved by applying the chained reduction, invoking a solver,
and using the chained reduction to retrieve a solution. For example,
if a user specifies a problem that is recognizably reducible to a quadratic
program, and if a quadratic solver is installed on the user's device, CVXPY 1.0
will automatically target it.

\section*{Acknowledgments}
We extend many thanks to Enzo Busseti, for significant contributions
to the quadratic program canonicalization procedure, and Yifan Lu, for fruitful
conversations about the analogy between our rewriting system and software
compilers.

\bibliographystyle{apacite}
\bibliography{cvxpy_rewriting}
\end{document}